\documentclass[12pt]{amsart}
\usepackage{amssymb}
\usepackage[all]{xy}
\usepackage{tikz}
\usepackage{amsmath}
\numberwithin{equation}{section}
\input xy
\xyoption{all}

\newtheorem{lemma}{Lemma}[section]

\newtheorem{theorem}[lemma]{Theorem}
\newtheorem{proposition}[lemma]{Proposition}

\theoremstyle{definition}
\newtheorem{remark}[lemma]{Remark}
\newtheorem{definition}[lemma]{Definition}

\newtheorem{example}[lemma]{Example}

\DeclareMathOperator{\Mod}{Mod}
\DeclareMathOperator{\modd}{mod}

\DeclareMathOperator{\Hom}{Hom}

\DeclareMathOperator{\Ext}{Ext}

\DeclareMathOperator{\Ker}{Ker}

\newtheorem{question}[lemma]{Question}
\newtheorem*{theorem a*}{Theorem A}
\newtheorem*{theorem b*}{Theorem B}

\begin{document}

\title{GABRIEL-QUILLEN EMBEDDING FOR $n$-EXACT CATEGORIES }

\author{Ramin Ebrahimi} 
\address{Department of Pure Mathematics\\
Faculty of Mathematics and Statistics\\
University of Isfahan\\
P.O. Box: 81746-73441, Isfahan, Iran\\ and School of Mathematics, Institute for Research in Fundamental Sciences (IPM), P.O. Box: 19395-5746, Tehran, Iran}
\email{ramin69@sci.ui.ac.ir / ramin.ebrahimi1369@gmail.com}



\begin{abstract}
Our first aim is to provide an analog of the Gabriel-Quillen embedding
theorem for $n$-exact categories. Also we give an example of an $n$-exact category that
in not an $n$-cluster tilting subcategory, and we suggest two possible ways for realizing
$n$-exact categories as $n$-cluster tilting subcategory. 
\end{abstract}

\maketitle

%%%%%%%%%%%%%%%%%%%%%%%%%%%%%%%%%%%%%%

\section{Introduction}
Higher Auslander-Reiten theory was introduced by Iyama in \cite{I1,I2}. It deals with $n$-cluster tilting subcategories of abelian and exact categories. Recently, Jasso \cite{J} introduced
$n$-abelian and $n$-exact categories as a higher-dimensional analogs of abelian and exact
categories, they are axiomatization of $n$-cluster tilting subcategories. Jasso proved that
each $n$-cluster tilting subcategory of an abelian (res, exact) category is $n$-abelian (res, $n$-exact).

In \cite{EN} and \cite{Kv}, independently it has been shown that any small $n$-abelian category is
equivalent to an $n$-cluster tilting subcategory of an abelian category. This note is an
attempt to generalize this result for $n$-exact categories. We give an example of an $n$-exact
category that is not equivalent to an $n$-cluster tilting subcategory, so we have to use a
different strategy for realizing $n$-exact categories as $n$-cluster tilting subcategories. 

Let $\mathcal{M}$ be a small $n$-exact category. We denote by $\rm Mod\mathcal{M}$ the category of all additive
contravariant functors from $\mathcal{M}$ to the category of all abelian groups. Let $\rm Eff(\mathcal{M})$ be
the subcategory of weakly effaceable functors, parallel to the proof of Gabriel-Quillen embedding
theorem we will show that composition of the Yoneda functor with the localisation functor 
\begin{equation}
\mathcal{M}\overset{Y}{\longrightarrow}\rm Mod\mathcal{M}\overset{q}{\longrightarrow}\frac{\rm Mod\mathcal{M}}{\rm Eff(\mathcal{M})} \notag
\end{equation}
sends $n$-exact sequences in $\mathcal{M}$ to exact sequences in $\mathcal{A}=\frac{\rm Mod\mathcal{M}}{\rm Eff(\mathcal{M})}$. Furthermore we will
show that this functor detects $n$-exact sequences and it's essential image is $n$-rigid in $\mathcal{A}$. In the end we suggest two possible ways for
realizing $n$-exact categories as $n$-cluster tilting subcategory. 

In section 2 we recall the definitions of $n$-exact categories, $n$-cluster tilting subcategories
and some of their basic properties. And we give an example of $n$-exact category that is
not an $n$-cluster tilting subcategory. In section 3 after recalling some results from localisation theory of
abelian categories, we construct the embedding $\mathcal{M}\hookrightarrow \mathcal{A}=\frac{\rm Mod\mathcal{M}}{\rm Eff(\mathcal{M})}$ with desired properties.
We end with a question that by results of this paper it make sense to has positive answer.

\subsection{Notation}
Throughout this paper, unless otherwise stated, $n$ always denotes a fixed positive integer and $\mathcal{M}$ is a fixed small $n$-exact category.

%%%%%%%%%%%%%%%%%%%%%%%%%%%%%%%%%%%%%%

\section{preliminaries}
In this section we recall the definition of $n$-exact category and $n$-cluster tilting subcategory. And we give an example of an $n$-exact category that can't be an $n$-cluster tilting subcategory of
an exact category 

\subsection{$n$-exact categories}
Let $\mathcal{M}$ be an additive category and $f:A\rightarrow B$ a morphism in $\mathcal{M}$. A weak cokernel of $f$ is a morphism $g:B\rightarrow C$ such that for all $C^{\prime} \in \mathcal{M}$  the sequence of abelian groups
\begin{equation}
\Hom(C,C')\overset{(g,C')}{\longrightarrow} \Hom(B,C')\overset{(f,C')}{\longrightarrow} \Hom(A,C') \notag
\end{equation}
is exact. The concept of a weak kernel is defined dually.

Let $d^0:X^0 \rightarrow X^1$ be a morphism in $\mathcal{M}$. An $n$-cokernel of $d^0$ is a sequence
\begin{equation}
(d^1, \ldots, d^n): X^1 \overset{d^1}{\rightarrow} X^2 \overset{d^2}{\rightarrow}\cdots \overset{d^{n-1}}{\rightarrow} X^n \overset{d^n}{\rightarrow} X^{n+1} \notag
\end{equation}
of objects and morphisms in $\mathcal{M}$ such that for each $Y\in \mathcal{M}$
the induced sequence of abelian groups
\begin{align}
0 \rightarrow \Hom(X^{n+1},Y) \rightarrow \Hom(X^n,Y) \rightarrow\cdots\rightarrow \Hom(X^1,Y) \rightarrow \Hom(X^0,Y) \notag
\end{align}
is exact. Equivalently, the sequence $(d^1, \ldots, d^n)$ is an $n$-cokernel of $d^0$ if for all $1\leq k\leq n-1$
the morphism $d^k$ is a weak cokernel of $d^{k-1}$, and $d^n$ is moreover a cokernel of $d^{n-1}$ \cite[Definition 2.2]{J}. The concept of an $n$-kernel of a morphism is defined dually.
\begin{definition}\label{d1}
Let $\mathcal{M}$ be an additive category. A left $n$-exact sequence in $\mathcal{M}$ is a complex
\begin{equation}
X^0 \overset{d^0}{\rightarrow} X^1 \overset{d^1}{\rightarrow} \cdots \overset{d^{n-1}}{\rightarrow} X^n \overset{d^n}{\rightarrow} X^{n+1} \notag
\end{equation}
such that $(d^0, \ldots, d^{n-1})$ is an $n$-kernel of $d^n$. The concept of right $n$-exact sequence is defined dually. An $n$-exact sequence is a sequence which is both a right $n$-exact sequence and a left $n$-exact sequence.
\end{definition}

Let
\begin{center} 
\begin{tikzpicture}
\node (X1) at (-4,1) {$X$};
\node (X2) at (-2,1) {$X^0$};
\node (X3) at (0,1) {$X^1$};
\node (X4) at (2,1) {$\ldots$};
\node (X5) at (4,1) {$X^{n-1}$};
\node (X6) at (6,1) {$X^n$};
\node (X7) at (-4,-1) {$Y$};
\node (X8) at (-2,-1) {$Y^0$};
\node (X9) at (0,-1) {$Y^1$};
\node (X10) at (2,-1) {$\ldots$};
\node (X11) at (4,-1) {$Y^{n-1}$};
\node (X12) at (6,-1) {$Y^n$};
\draw [->,thick] (X1) -- (X7) node [midway,left] {$f$};
\draw [->,thick] (X2) -- (X8) node [midway,left] {$f^0$};
\draw [->,thick] (X3) -- (X9) node [midway,left] {$f^1$};
\draw [->,thick] (X5) -- (X11) node [midway,left] {$f^{n-1}$};
\draw [->,thick] (X6) -- (X12) node [midway,left] {$f^n$};
\draw [->,thick] (X2) -- (X3) node [midway,above] {$d_X^0$};
\draw [->,thick] (X3) -- (X4) node [midway,above] {$d_X^1$};
\draw [->,thick] (X4) -- (X5) node [midway,above] {$d_X^{n-2}$};
\draw [->,thick] (X5) -- (X6) node [midway,above] {$d_X^{n-1}$};
\draw [->,thick] (X8) -- (X9) node [midway,above] {$d_Y^0$};
\draw [->,thick] (X9) -- (X10) node [midway,above] {$d_Y^1$};
\draw [->,thick] (X10) -- (X11) node [midway,above] {$d_Y^{n-2}$};
\draw [->,thick] (X11) -- (X12) node [midway,above] {$d_Y^{n-1}$};
\end{tikzpicture}
\end{center}
be a morphism of complexes in an additive category. The mapping cone
$C = C( f )$
is the complex
\begin{equation} \label{Cone}
X^0 \xrightarrow{d_C^{-1}} X^1\oplus Y^0 \xrightarrow{d_C^0} \ldots \xrightarrow{d_C^{n-2}} X^n\oplus Y^{n-1}\xrightarrow{d_C^{n-1}} Y^n,
\end{equation}
where
\begin{center}
$d_C^k:=\begin{pmatrix}
-d_X^{k+1}& 0\\\\
f^{k+1} & d_Y^k
\end{pmatrix}: X^{k+1}\oplus Y^k\to X^{k+2} \oplus Y^{k+1}$
\end{center}
for each
$k\in \{-1, 0, \ldots, n-1\}$.
In particular
$d_C^{-1}=\begin{pmatrix}
-d_X^0\\
f^0
\end{pmatrix}$
and
$d_C^{n-1}=(f^n d_Y^{n-1})$.
\begin{itemize}
\item
We say that the above diagram is an $n$-pull back of
$Y$
along
$f^n$
if \eqref{Cone} is a left $n$-exact sequence.
\item
We say that the above diagram is an $n$-push out of
$X$
along
$f^0$
if \eqref{Cone} is a right $n$-exact sequence.
\end{itemize}

Let $\mathcal{A}$ be an additive category and $\mathcal{B}$ be a full subcategory of $\mathcal{A}$. $\mathcal{B}$ is called
covariantly finite in $\mathcal{A}$ if for every $A\in \mathcal{A}$ there exists an object $B\in\mathcal{B}$ and a morphism
$f : A\rightarrow B$ such that, for all $B'\in\mathcal{B}$, the sequence of abelian groups $\Hom_{\mathcal{A}}(B, B')\rightarrow \Hom_{\mathcal{A}}(A, B')\rightarrow 0$ is exact. Such a morphism $f$ is called a left $\mathcal{B}$-approximation of $A$. The notions of contravariantly
finite subcategory of $\mathcal{A}$ and right $\mathcal{B}$-approximation are defined dually. A functorially
finite subcategory of $\mathcal{A}$ is a subcategory which is both covariantly and contravariantly finite
in $\mathcal{A}$.

Let $X$ and $Y$ be two $n$-exact sequences. We remained that a morphism $f:X\rightarrow Y$ of $n$-exact sequences is a morphism of complexes.
We say that a morphism $f:X\rightarrow Y$ of $n$-exact sequences is a weak isomorphism if $f^k$
and $f^{k+1}$ are isomorphisms for some $k\in \{0,1, . . . , n + 1\}$ with $n + 2 := 0$.

\begin{definition}$($\cite[Definition 4.2]{J}$)$
Let $\mathcal{M}$ be an additive category. An $n$-exact structure
on $\mathcal{M}$ is a class $\mathcal{X}$ of $n$-exact sequences in $\mathcal{M}$, closed under weak isomorphisms of $n$-exact
sequences, and which satisfies the following axioms:
\begin{itemize}
\item[$(E0)$]
The sequence $0\rightarrowtail 0\rightarrow \cdots\rightarrow 0\twoheadrightarrow 0$ is an $\mathcal{X}$-admissible $n$-exact sequence. 
\item[$(E1)$]
The class of $\mathcal{X}$-admissible monomorphisms is closed under composition. 
\item[$(E1^{op})$]
The class of $\mathcal{X}$-admissible epimorphisms is closed under composition. 
\item[$(E2)$]
For each $\mathcal{X}$-admissible $n$-exact sequence $X$ and each morphism $f:X^0\rightarrow Y^0$, there exists an $n$-pushout diagram of $(d_X^0,\cdots , d_X^{n-1})$ along $f$ such that $d_Y^0$ is an $\mathcal{X}$-admissible monomorphism. The situation is illustrated in the following commutative diagram:
\begin{center}
\begin{tikzpicture}
\node (X1) at (-4,1) {$X^0$};
\node (X2) at (-2,1) {$X^1$};
\node (X3) at (0,1) {$\cdots$};
\node (X4) at (2,1) {$X^n$};
\node (X5) at (4,1) {$X^{n+1}$};
\node (X6) at (-4,-1) {$Y^0$};
\node (X7) at (-2,-1) {$Y^1$};
\node (X8) at (0,-1) {$\cdots$};
\node (X9) at (2,-1) {$X^n$};
\draw [>->,thick] (X1) -- (X2) node [midway,above] {$d_X^0$};
\draw [->,thick] (X2) -- (X3) node [midway,above] {$d_X^1$};
\draw [->,thick] (X3) -- (X4) node [midway,above] {$d_X^{n-1}$};
\draw [->>,thick] (X4) -- (X5) node [midway,above] {$d_X^n$};
\draw [>->,thick,dashed] (X6) -- (X7) node [midway,above] {$d_Y^0$};
\draw [->,thick,dashed] (X7) -- (X8) node [midway,above] {$d_Y^1$};
\draw [->,thick,dashed] (X8) -- (X9) node [midway,above]{$d_Y^{n-1}$};
\draw [->,thick] (X1) -- (X6) node [midway,left] {$f$};
\draw [->,thick,dashed] (X2) -- (X7) node [midway,left] {};
\draw [->,thick,dashed] (X4) -- (X9) node [midway,left] {};
\end{tikzpicture}
\end{center}

\item[$(E2^{op})$]
For each $\mathcal{X}$-admissible $n$-exact sequence $Y$ and each morphism $g:X^{n+1}\rightarrow Y^{n+1}$, there exists an $n$-pull back diagram of $(d_Y^1,\cdots , d_Y^{n})$ along $g$ such that $d_X^{n}$ is an $\mathcal{X}$-admissible epimorphism. The situation is illustrated in the following commutative diagram:
\begin{center}
\begin{tikzpicture}
%\node (X1) at (-4,1) {$X^0$};
\node (X2) at (-2,1) {$X^1$};
\node (X3) at (0,1) {$\cdots$};
\node (X4) at (2,1) {$X^n$};
\node (X5) at (4,1) {$X^{n+1}$};
\node (X6) at (-4,-1) {$Y^0$};
\node (X7) at (-2,-1) {$Y^1$};
\node (X8) at (0,-1) {$\cdots$};
\node (X9) at (2,-1) {$Y^n$};
\node (X10) at (4,-1) {$Y^{n+1}$};
%\draw [->,thick] (X1) -- (X2) node [midway,above] {$d_X^0$};
\draw [->,thick,dashed] (X2) -- (X3) node [midway,above] {$d_X^1$};
\draw [->,thick,dashed] (X3) -- (X4) node [midway,above] {$d_X^{n-1}$};
\draw [->>,thick,dashed] (X4) -- (X5) node [midway,above] {$d_X^n$};
\draw [>->,thick] (X6) -- (X7) node [midway,above] {$d_Y^0$};
\draw [->,thick] (X7) -- (X8) node [midway,above] {$d_Y^1$};
\draw [->,thick] (X8) -- (X9) node [midway,above] {$d_Y^{n-1}$};
\draw [->>,thick] (X9) -- (X10) node [midway,above] {$d_Y^{n}$};
\draw [->,thick] (X5) -- (X10) node [midway,right] {$g$};
\draw [->,thick,dashed] (X2) -- (X7) node [midway,left] {};
\draw [->,thick,dashed] (X4) -- (X9) node [midway,left] {};
\end{tikzpicture}
\end{center}
\end{itemize}
\end{definition}

An $n$-exact category is a pair $(\mathcal{M},\mathcal{X})$ where $\mathcal{M}$ is an additive category and $\mathcal{X}$ is an
$n$-exact structure on $\mathcal{M}$. If the class $\mathcal{X}$ is clear from the context, we identify $\mathcal{M}$ with the
pair $(\mathcal{M},\mathcal{X})$. The members of $\mathcal{X}$ are called $\mathcal{X}$-admissible $n$-exact sequences, or simply
admissible $n$-exact sequences when $\mathcal{X}$ is clear from the context. Furthermore, if 
\begin{equation}
X^0 \overset{d^0}{\rightarrowtail} X^1 \overset{d^1}{\rightarrow} \cdots \overset{d^{n-1}}{\rightarrow} X^n \overset{d^n}{\twoheadrightarrow} X^{n+1} \notag
\end{equation}
is an admissible $n$-exact sequence, $d^0$ is called admissible monomorphism and $d^n$ is called
admissible epimorphism. 

\begin{definition}$($\cite[Definition 4.13]{J}$)$
Let $(\mathcal{E},\mathcal{X})$ be an exact category and $\mathcal{M}$ a subcategory of $\mathcal{E}$. $\mathcal{M}$ is called an $n$-cluster tilting subcategory of $(\mathcal{E},\mathcal{X})$ if the following conditions are satisfied.
\begin{itemize}
\item[$(i)$]
Every object $E\in \mathcal{E}$ has a left $\mathcal{M}$-approximation by an $\mathcal{X}$-admissible monomorphism $E\rightarrowtail M$.
\item[$(ii)$]
Every object $E\in \mathcal{E}$ has a right $\mathcal{M}$-approximation by an $\mathcal{X}$-admissible epiomorphism $M'\twoheadrightarrow E$.
\item[$(iii)$]
We have
\begin{align}
\mathcal{M}& = \{ E\in \mathcal{E} \mid \forall i\in \{1, \ldots, n-1 \}, \Ext_{\mathcal{E}}^i(E,\mathcal{M})=0 \}\notag \\
& =\{ E\in \mathcal{E} \mid \forall i\in \{1, \ldots, n-1 \}, \Ext_{\mathcal{E}}^i(\mathcal{M},E)=0 \}.\notag
\end{align}
\end{itemize}
Note that $\mathcal{E}$ itself is the unique $1$-cluster tilting subcategory of $\mathcal{E}$.
\end{definition}

A full subcategory $\mathcal{M}$ of an exact or abelian category $\mathcal{E}$ is called $n$-rigid, if for every
two objects $M,N\in \mathcal{M}$ and for every $k\in \{1,\cdots, n - 1\}$, we have $\Ext_{\mathcal{E}}^i(\mathcal{M},\mathcal{M})=0$. Any
$n$-cluster tilting subcategory $\mathcal{M}$ of an exact category $\mathcal{E}$ is $n$-rigid. 

The following theorem gives the main source of $n$-exact categories. 

\begin{theorem}$($\cite[Theorem 4.14]{J}$)$\label{Tne}
Let $(\mathcal{E},\mathcal{X})$ be an exact category and $\mathcal{M}$ be an $n$-cluster
tilting subcategory of $(\mathcal{E},\mathcal{X})$. Let $\mathcal{Y}=\mathcal{Y}(\mathcal{M},\mathcal{X})$ be the class of all $\mathcal{X}$-acyclic complexes 
\begin{equation}
X^0 \overset{d^0}{\rightarrowtail} X^1 \overset{d^1}{\rightarrow} \cdots \overset{d^{n-1}}{\rightarrow} X^n \overset{d^n}{\twoheadrightarrow} X^{n+1} \notag
\end{equation}
such that for all $k\in\{0, 1,\cdots, n + 1\}$ we have $X^k\in \mathcal{M}$. Then $(\mathcal{M},\mathcal{Y})$ is an $n$-exact category. 
\end{theorem}

Let $\mathcal{M}$ be an additive category and $M$ be an object of $\mathcal{M}$. A morphism $e\in \mathcal{M}(M,M)$
is called idempotent if $e^2=e$. $\mathcal{M}$ is called idempotent complete if for every idempotent $e\in
\mathcal{M}(M, M)$ there exist an object $N$ and morphisms $f\in \mathcal{M}(M, N)$ and $g\in \mathcal{M}(N, M)$
such that $gf = e$ and $fg = 1_N$. Assume that $r : M \rightarrow M'$ is a retraction with section
$s : M' \rightarrow M$. Then $sr : M \rightarrow M$ is an idempotent. It is well known that if $r : M \rightarrow M'$
has a kernel $k : K \rightarrow M$, this idempotent splits and there is a canonical isomorphism
$M\cong K\oplus M'$ \cite{Bu}.

In abelian categories all retractions have kernels, but in exact categories this does not happen 
in general. An exact category where
all retractions have kernels are called weakly
idempotent complete \cite{Bu}. But it is obvious that any admissible epimorphism in an exact category, that is a
retraction has a kernel. 

Let $\mathcal{M}$ be an $n$-cluster tilting subcategory of an exact category $(\mathcal{E},\mathcal{X})$, and  $\mathcal{Y}$ be the class of all $\mathcal{X}$-acyclic complexes
\begin{equation}
X^0 \overset{d^0}{\rightarrowtail} X^1 \overset{d^1}{\rightarrow} \cdots \overset{d^{n-1}}{\rightarrow} X^n \overset{d^n}{\twoheadrightarrow} X^{n+1} \notag
\end{equation}
such that for all $k\in\{0, 1,\cdots, n + 1\}$ we have $X^k\in \mathcal{M}$. By Theorem \ref{Tne}, $(\mathcal{M},\mathcal{Y})$ is an $n$-exact category. If $M,N\in \mathcal{M}$, a morphism
$f : M \rightarrow N$ is $\mathcal{Y}$-admissible epimprphism if and only if it is $\mathcal{X}$-admissible epimorphism \cite{J}. Thus if $f : M \rightarrow N$ is an $\mathcal{Y}$-admissible epimorphism that is a retraction with section $g : N \rightarrow M$, the idempotent $gf : M \rightarrow M$ splits and $M\cong N\oplus \rm Ker(f)$. By the definition of $n$-cluster tilting subcategory, $\rm Ker(f) \in \mathcal{M}$. 

\begin{example}$($\cite[Example 3.5]{J}$)$\label{Ex}
Let $n \geq 2$ and $K$ be a field. Consider the full subcategory $\mathbb{V}$ of $\modd K$ given by the finite dimensional $K$-vector spaces of dimension different
from $1$. Then it has been shown in \cite[Example 3.5]{J} that $\mathbb{V}$ is not idempotent complete,
but it satisfies other axioms of $n$-abelian category. By a similar argument the class of
all exact sequences with $n + 2$ term is an $n$-exact structure on $\mathbb{V}$. But there exist an
admissible epimorphism $K^3 \rightarrow K^2$ which is a retraction, that doesn't give a splitting of
$K^3$. Thus $\mathbb{V}$ can't be an $n$-cluster tilting subcategory. Note that we can consider $\mathbb{V}$ as an
$n$-cluster tilting subcategory of itself, but in this case the induced $n$-exact structure is
different than the class of all exact sequences in $\modd K$. 
\end{example}

%%%%%%%%%%%%%%%%%%%%%%%%%%%%%%%%%%%%%%

\section{Embeddings into abelian categories}
Let $\mathcal{M}$ be a small $n$-exact category. In this section we find an abelian category $\mathcal{A}$
and an embedding $H : \mathcal{M} \hookrightarrow \mathcal{A}$, such that $H$ sends $n$-exact sequences in $\mathcal{M}$ to exact
sequences in $\mathcal{A}$. Furthermore we will show that $H$ detects $n$-exact sequences and it's essential image is $n$-rigid in $\mathcal{A}$. 

First we recall localisation theory of abelian categories, for reader can find proof in standard
textbooks or Gabriel thesis \cite{G}.
Let $\mathcal{A}$ be an abelian category. A subcategory $\mathcal{C}$ of $\mathcal{A}$ is called a \textbf{Serre subcategory} if for
any exact sequence
\begin{equation}
0\rightarrow A_1\rightarrow A_2\rightarrow A_3\rightarrow 0 \notag
\end{equation}
we have that $A_2\in \mathcal{C}$ if and only if $A_1\in \mathcal{C}$ and $A_3\in \mathcal{C}$. In this case we have the quotient category $\dfrac{\mathcal{A}}{\mathcal{C}}$ that is by definition localisation of $\mathcal{A}$ with respect to the class of all morphisms $f:X\rightarrow Y$ such that $\rm Ker(f), \rm Coker(f)\in \mathcal{C}$.

\begin{theorem}\label{T1}
Let $\mathcal{C}$ be a Serre subcategory of $\mathcal{A}$, and let $q:\mathcal{A}\rightarrow \dfrac{\mathcal{A}}{\mathcal{C}}$ denote the canonical
functor to the localization. The following
statements hold:
\begin{itemize}
\item[$(i)$]
$\dfrac{\mathcal{A}}{\mathcal{C}}$ is an abelian category and $q$ is an exact functor.
\item[$(ii)$]
$q(C) = 0$ for all $C\in \mathcal{C}$, and any exact functor $F:\mathcal{A}\rightarrow \mathcal{D}$ annihilating $\mathcal{C}$ where $\mathcal{D}$ is abelian must factor uniquely through $q$. 
\end{itemize}
\end{theorem}

A Serre subcategory $\mathcal{C}\subseteq \mathcal{A}$ is called a \textbf{localizing subcategory} if the canonical functor $q:\mathcal{A}\rightarrow \dfrac{\mathcal{A}}{\mathcal{C}}$ admits a right adjoint $r:\dfrac{\mathcal{A}}{\mathcal{C}}\rightarrow \mathcal{A}$. The right adjoint $r$ is called the \textbf{section functor}, which always is fully faithful. Note that a localising subcategory is closed under all coproducts which exist in $\mathcal{A}$. The converse is true for Grothendieck categories, indeed we have the following result. 

\begin{theorem}\label{T2}
Let $\mathcal{C}$ be a Serre subcategory of a Grothendieck category $\mathcal{A}$. The following
statements hold:
\begin{itemize}
\item[$(i)$]
$\mathcal{C}$ is a localising subcategory if and only if it is closed under coproducts. 
\item[$(ii)$]
In this case the quotient category $\dfrac{\mathcal{A}}{\mathcal{C}}$ is a Grothendieck category. 
\end{itemize}
\end{theorem}

Let $\mathcal{C}$ be a Serre subcategory of an abelian category $\mathcal{A}$. Recall that an object $A\in \mathcal{A}$ is called $\mathcal{C}$-\textbf{closed}
if for every morphism $f : X \rightarrow Y$ with $\Ker(f)\in \mathcal{C}$ and $\rm Coker(f) \in \mathcal{C}$ we have that
$\Hom_{\mathcal{A}}(f, A)$ is bijective. Denote by $\mathcal{C}^{\bot}$ the full subcategory of all $\mathcal{C}$-closed objects, the following result is well known. 

\begin{theorem}\label{T3}
Let $\mathcal{C}$ be a Serre subcategory of an abelian category $\mathcal{A}$. The following
statements hold:
\begin{itemize}
\item[$(i)$]
We have
\begin{center}
$\mathcal{C}^{\bot}=\{A\in \mathcal{A} \mid \Hom(\mathcal{C},A)=0=\Ext^1(\mathcal{C},A)\}.$
\end{center}
\item[$(ii)$]
For $A\in \mathcal{A}$ and $B\in \mathcal{C}^{\bot}$, the natural homomorphism $q_{A,B}:\Hom_{\mathcal{A}}(A,B)\rightarrow \Hom_{\frac{\mathcal{A}}{\mathcal{C}}}(q(A),q(B))$ is an isomorpism.
\item[$(iii)$]
If $\mathcal{C}$ is a localizing subcategory, the restriction $q:\mathcal{C}^{\bot}\rightarrow \dfrac{\mathcal{A}}{\mathcal{C}}$ is an equivalence of categories.
\item[$(iv)$]
If $\mathcal{C}$ is localising and $\mathcal{A}$ has injective envelopes, then $\mathcal{C}^{\bot}$ has injective envelopes and the inclusion functor $\mathcal{C}^{\bot}\hookrightarrow \mathcal{A}$ preserves injective envelopes. 
\end{itemize}
\end{theorem}

We also need the following technical lemma. 

\begin{lemma}\label{L4}
Let $0\rightarrow A\rightarrow L\rightarrow M\rightarrow 0$ be an exact sequence in $\mathcal{A}$ with $L\in \mathcal{C}^{\bot}$, then
$A\in \mathcal{C}^{\bot}$ if and only if $\Hom(\mathcal{C}, M) = 0$. 
\end{lemma}

Now we want to apply the above general results to $\rm{Mod}\mathcal{M}$, where $\mathcal{M}$ is a small $n$-exact category. Recall that $\rm{Mod}\mathcal{M}$ is
the category of all additive contravariant functors from $\mathcal{M}$ to the category of all abelian
groups. It is an abelian category with all limits and colimits, which are defined point-wise . Also by the Yoneda's lemma, representable functors are projective and the direct sum of all representable functors $\Sigma_{X\in \mathcal{M}}\Hom(-,X)$, is a generator for $\rm{Mod}\mathcal{M}$. Thus $\rm{Mod}\mathcal{M}$ is a Grothendieck category.

A functor $F\in \rm{Mod}\mathcal{M}$ is called \textbf{weakly effaceable}, if for each object $X\in \mathcal{M}$ and $x\in F(X)$
there exists an admissible epimorphism $f : Y \rightarrow X$ such that $F(f)(x) = 0$. We denote by
$\rm Eff(\mathcal{M})$ the full subcategory of all weakly effaceable functors. For each $k\in \{1,\cdots, n\}$ we denote by
$\mathcal{L}_k(\mathcal{M})$ the full subcategory of $\rm{Mod}\mathcal{M}$ consist of all functors like $F$ such that for every $n$-exact sequence
\begin{equation}
X^0 \overset{}{\rightarrowtail} X^1 \overset{}{\rightarrow} \cdots \overset{}{\rightarrow} X^n \overset{}{\twoheadrightarrow} X^{n+1} \notag
\end{equation}
the sequence of abelian groups
\begin{equation}
0 \overset{}{\rightarrow} F(X^{n+1}) \overset{}{\rightarrow}F(X^n)\rightarrow \cdots \overset{}{\rightarrow} F(X^{n-k})\notag
\end{equation}
is exact. Also for a Serre subcategory $\mathcal{C}$ of an abelian category $\mathcal{A}$ we set $\mathcal{C}^{{\bot}_k} = \{A \in \mathcal{A} | \Ext^{0,...,k}(\mathcal{C}, A) = 0\}$. Note that $\mathcal{C}^{{\bot}_1}=\mathcal{C}^{\bot}$ by Theorem \ref{T3}.

\begin{proposition}\label{P5}
\begin{itemize}
\item[(i)]
$\rm Eff(\mathcal{M})$ is a localizing subcategory of $\rm{Mod}\mathcal{M}$.
\item[(ii)]
$\rm Eff(\mathcal{M})^{\bot}=\mathcal{L}_1(\mathcal{M})$.
\end{itemize}
\begin{proof}
$(i)$ We need to show that $\rm Eff(\mathcal{M})$
is a Serre subcategory closed under coproducts, because $\rm{Mod}\mathcal{M}$ is a Grothendieck category. The proof is similar to the classical case of exact categories. We only prove that $\rm Eff(\mathcal{M})$ is closed under extensions. Let
\begin{center}
$0\rightarrow F_1\overset{\alpha}{\rightarrow}F_2\overset{\beta}{\rightarrow} F_3\rightarrow 0$
\end{center}
 be a short exact sequence in $\rm{Mod}\mathcal{M}$
and
$F_1,F_3\in \rm{Eff}(\mathcal{M})$.
We want to show that
$F_2\in \rm{Eff}(\mathcal{M})$.
Let
$X\in \mathcal{M}$
and
$x_2\in F_2(X)$.
Set
$x_3=\beta_X(x_2)\in F_3(X)$.
By assumption there exist an admissible epimorphism
$f:Y\rightarrow X$
such that
$F_3(f)(x_3)=0$.
\begin{center}
\begin{tikzpicture}
\node (X1) at (-6,1) {$0$};
\node (X2) at (-3,1) {$F_1(X)$};
\node (X3) at (0,1) {$F_2(X)$};
\node (X4) at (3,1) {$F_3(X)$};
\node (X5) at (6,1) {$0$};
\node (X6) at (-6,-1) {$0$};
\node (X7) at (-3,-1) {$F_1(Y)$};
\node (X8) at (0,-1) {$F_2(Y)$};
\node (X9) at (3,-1) {$F_3(Y)$};
\node (X10) at (6,-1) {$0$};
\draw [->,thick] (X1) -- (X2) node [midway,above] {};
\draw [->,thick] (X2) -- (X3) node [midway,above] {$\alpha_X$};
\draw [->,thick] (X3) -- (X4) node [midway,above] {$\beta_X$};
\draw [->,thick] (X4) -- (X5) node [midway,above] {};
\draw [->,thick] (X2) -- (X7) node [midway,left] {$F_1(f)$};
\draw [->,thick] (X3) -- (X8) node [midway,left] {$F_2(f)$};
\draw [->,thick] (X4) -- (X9) node [midway,left] {$F_3(f)$};
\draw [->,thick] (X6) -- (X7) node [midway,left] {};
\draw [->,thick] (X7) -- (X8) node [midway,above] {$\alpha_Y$};
\draw [->,thick] (X8) -- (X9) node [midway,above] {$\beta_Y$};
\draw [->,thick] (X9) -- (X10) node [midway,left] {};
\end{tikzpicture}
\end{center}
Using the above commutative diagram
$F_2(f)(x_2)\in \rm{Ker}(\beta_Y)=\rm{Im}(\alpha_Y)$.
Thus there exists
$y_1\in F_1(Y)$
such that
$\alpha_Y(y_1)=F_2(f)(x_2)$.
Again by assumption there exist an admissible epimorphism
$g:Z\rightarrow Y$
such that
$F_1(g)(y_1)=0$.
\begin{center}
\begin{tikzpicture}
\node (X1) at (-6,1) {$0$};
\node (X2) at (-3,1) {$F_1(Y)$};
\node (X3) at (0,1) {$F_2(Y)$};
\node (X4) at (3,1) {$F_3(Y)$};
\node (X5) at (6,1) {$0$};
\node (X6) at (-6,-1) {$0$};
\node (X7) at (-3,-1) {$F_1(Z)$};
\node (X8) at (0,-1) {$F_2(Z)$};
\node (X9) at (3,-1) {$F_3(Z)$};
\node (X10) at (6,-1) {$0$};
\draw [->,thick] (X1) -- (X2) node [midway,above] {};
\draw [->,thick] (X2) -- (X3) node [midway,above] {$\alpha_Y$};
\draw [->,thick] (X3) -- (X4) node [midway,above] {$\beta_Y$};
\draw [->,thick] (X4) -- (X5) node [midway,above] {};
\draw [->,thick] (X2) -- (X7) node [midway,left] {$F_1(g)$};
\draw [->,thick] (X3) -- (X8) node [midway,left] {$F_2(g)$};
\draw [->,thick] (X4) -- (X9) node [midway,left] {$F_3(g)$};
\draw [->,thick] (X6) -- (X7) node [midway,left] {};
\draw [->,thick] (X7) -- (X8) node [midway,above] {$\alpha_Z$};
\draw [->,thick] (X8) -- (X9) node [midway,above] {$\beta_Z$};
\draw [->,thick] (X9) -- (X10) node [midway,left] {};
\end{tikzpicture}
\end{center}
Using the above commutative diagram
$F_2(gf)(x_2)=F_2(g)F_2(f)(x_2)=F_2(g)\alpha_Y(y_1)=\alpha_ZF_1(g)(y_1)=0$. Since $gf$ is an admissible epimorphism, $F_2\in \rm Eff(\mathcal{M})$

$(ii)$ Let $L\in \mathcal{L}_1(\mathcal{M})$, consider the exact sequence $0 \rightarrow L \rightarrow I \rightarrow M \rightarrow 0$ where I is injective envelope of $L$. First note that for every $n$-exact sequence $X^0\rightarrowtail X^1\rightarrow \cdots\rightarrow X^n\twoheadrightarrow X^{n+1}$ by definition
\begin{center}
$0 \rightarrow (-,X^0) \rightarrow (-,X^1) \rightarrow \cdots \rightarrow (-,X^n)
\rightarrow (-,X^{n+1}) $
\end{center}
is exact, applying the exact functor $(-, I)$ to this sequence we obtain that $ I(X^{n+1}) \rightarrow
I(X^n) \rightarrow\cdots\rightarrow I(X^1) \rightarrow I(X^0)\rightarrow 0$ is exact. Also because $L\in \mathcal{L}_1(\mathcal{M})$ it doesn't have any nonzero weakly effaceable
subobject, so $I$ doesn't have any nonzero weakly effaceable subobject because it is an injective envelope if $L$. This
means that $I$ is an $n$-exact functor i.e
\begin{equation}
0 \rightarrow I(X^{n+1}) \rightarrow I(X^n) \rightarrow\cdots\rightarrow I(X^1) \rightarrow I(X^0) \rightarrow 0 \notag 
\end{equation}
is exact for all $n$-exact sequences in $\mathcal{M}$. Consider the following commutative diagram.
\begin{center}
\begin{tikzpicture}
\node (X1) at (-3,1.5) {$0$};
\node (X2) at (0,1.5) {$0$};
\node (X3) at (3,1.5) {$0$};
\node (X4) at (-6,0) {$0$};
\node (X5) at (-3,0) {$L(X^{n+1})$};
\node (X6) at (0,0) {$I(X^{n+1})$};
\node (X7) at (3,0) {$M(X^{n+1})$};
\node (X8) at (6,0) {$0$};
\node (X9) at (-6,-1.5) {$0$};
\node (X10) at (-3,-1.5) {$L(X^{n})$};
\node (X11) at (0,-1.5) {$I(X^{n})$};
\node (X12) at (3,-1.5) {$M(X^{n})$};
\node (X13) at (6,-1.5) {$0$};
\node (X14) at (-6,-3) {$0$};
\node (X15) at (-3,-3) {$L(X^{n-1})$};
\node (X16) at (0,-3) {$I(X^{n-1})$};
\node (X17) at (3,-3) {$M(X^{n-1})$};
\node (X18) at (6,-3) {$0$};
\draw [->,thick] (X4) -- (X5) node [midway,left] {};
\draw [->,thick] (X5) -- (X6) node [midway,left] {};
\draw [->,thick] (X6) -- (X7) node [midway,left] {};
\draw [->,thick] (X7) -- (X8) node [midway,left] {};
\draw [->,thick] (X9) -- (X10) node [midway,left] {};
\draw [->,thick] (X10) -- (X11) node [midway,left] {};
\draw [->,thick] (X11) -- (X12) node [midway,above] {};
\draw [->,thick] (X12) -- (X13) node [midway,left] {};
\draw [->,thick] (X14) -- (X15) node [midway,left] {};
\draw [->,thick] (X15) -- (X16) node [midway,left] {};
\draw [->,thick] (X16) -- (X17) node [midway,left] {};
\draw [->,thick] (X17) -- (X18) node [midway,left] {};
\draw [->,thick] (X1) -- (X5) node [midway,left] {};
\draw [->,thick] (X5) -- (X10) node [midway,above] {};
\draw [->,thick] (X10) -- (X15) node [midway,left] {};
\draw [->,thick] (X2) -- (X6) node [midway,left] {};
\draw [->,thick] (X6) -- (X11) node [midway,left] {};
\draw [->,thick] (X11) -- (X16) node [midway,left] {};
\draw [->,thick] (X3) -- (X7) node [midway,left] {};
\draw [->,thick] (X7) -- (X12) node [midway,left] {};
\draw [->,thick] (X12) -- (X17) node [midway,above] {};
\end{tikzpicture}
\end{center}
All rows are exact by assumption, and the left-hand and middle columns are exact,
now long exact sequence theorem \cite[Theorem 1.3.1]{We} tells that $0 \rightarrow M(X^{n+1}) \rightarrow M(X^n)$ is exact. Thus
$\rm Hom(Eff(\mathcal{M}),M) = 0$. Now by Lemma \ref{L4} $L \in \rm Eff(\mathcal{M})^{\bot}$. For the converse inclusion $\rm Eff(\mathcal{M})^{\bot}\subseteq \mathcal{L}_1(\mathcal{M})$, let $L\in \rm Eff(\mathcal{M})^{\bot}$ and consider the short exact sequence $0\rightarrow L\rightarrow I\rightarrow M\rightarrow 0$ where $I$ is an injective envelope of $L$. Thus by Lemma \ref{L4} $\rm Hom(Eff(\mathcal{M}),M) = 0$, that means $0 \rightarrow M(X^{n+1}) \rightarrow M(X^n)$ is exact. Again by long exact sequence theorem, the left-hand column is exact.
\end{proof}
\end{proposition}

The following observation is interesting and is our motivation for Question \ref{Q}.

\begin{proposition}\label{P6}
For every $k\in \{1,\cdots, n\}$, $\rm Eff(\mathcal{M})^{\bot_k}=\mathcal{L}_k(\mathcal{M})$. 
\begin{proof}
We want to prove by induction that for all $1\leq k\leq n$, $\rm Eff(\mathcal{M})^{\bot_k}=\mathcal{L}_k(\mathcal{M})$. By Proposition \ref{P5} $\rm Eff(\mathcal{M})^{\bot_1}=\mathcal{L}_1(\mathcal{M})$. Let $k\geq 2$, $L\in \rm Eff(\mathcal{M})^{\bot_1}=\mathcal{L}_1(\mathcal{M})$ and
\begin{center}
$X:X^0\rightarrowtail X^1\rightarrow \cdots\rightarrow X^n\twoheadrightarrow X^{n+1}$
\end{center}
be an arbitrary $n$-exact sequence. Consider the exact sequence $0 \rightarrow L \rightarrow I \rightarrow M \rightarrow 0$ where $I$ is injective envelope of $L$. Note that as we see in the proof of Proposition \ref{P5}
\begin{equation}
0 \rightarrow I(X^{n+1}) \rightarrow I(X^n) \rightarrow\cdots\rightarrow I(X^1) \rightarrow I(X^0) \rightarrow 0 \notag 
\end{equation}
is exact. By dimension shifting $L\in \rm Eff(\mathcal{M})^{\bot_k}$ if and only if $M\in \rm Eff(\mathcal{M})^{\bot_{k-1}}=\mathcal{L}_{k-1}(\mathcal{M})$. Applying the long exact sequence theorem \cite[Theorem 1.3.1]{We} to the following
short exact sequence of complexes.
\begin{center}
$0 \rightarrow L(X) \rightarrow I(X) \rightarrow M(X) \rightarrow 0 $
\end{center}
Because the middle column is exact we obtain that $M\in \mathcal{L}_{k-1}(\mathcal{M})$ if and only if $L\in \mathcal{L}_k(\mathcal{M})$.
\end{proof}
\end{proposition}

We denote by $H:\mathcal{M}\rightarrow \mathcal{L}_1(\mathcal{M})$ the composition of the Yoneda functor $\mathcal{M}\rightarrow \Mod\mathcal{M}$ with the localization functor $\Mod \mathcal{M}\rightarrow \frac{\Mod \mathcal{M}}{\rm{Eff}(\mathcal{M})}\simeq \rm{Eff}(\mathcal{M})^{\bot}=\mathcal{L}_1(\mathcal{M})$. Thus $H(X)=(-,X):\mathcal{M}^{op}\rightarrow \rm Ab$. For simplicity we denote $(-,X)$ by $H_X$.

\begin{proposition}\label{PM}
\begin{itemize}
\item[(i)]
For every $n$-exact sequence $X^0\rightarrowtail X^1\rightarrow \cdots\rightarrow X^n\twoheadrightarrow X^{n+1}$ in $\mathcal{M}$,
\begin{center}
$0\rightarrow H_{X^0}\rightarrow H_{X^1}\rightarrow \cdots\rightarrow H_{X^{n}}\rightarrow H_{X^{n+1}}\rightarrow 0$
\end{center}
is exact in $\mathcal{L}_1(\mathcal{M})$.

\item[(ii)]
The essential image of $H:\mathcal{M}\rightarrow \mathcal{L}_1(\mathcal{M})$ is $n$-rigid.
\end{itemize}
\begin{proof}
Because the cokernel of $H_{X^n}\rightarrow H_{X^{n+1}}$ is weakly effaceable, (i) follows.

Let $X,Y\in \mathcal{M}$ and $H_Y\rightarrow I^0$ be the injective envelope of $H_Y$ in $\Mod \mathcal{M}$.
Because $H_Y=(-,Y)\in \mathcal{L}_n(\mathcal{M})$, by the proofs of Proposition \ref{P5} and Proposition \ref{P6} $I^0\in \mathcal{L}_n(\mathcal{M})$ and in the short exact sequence
\begin{center}
$0\rightarrow H_Y\rightarrow I^0\rightarrow \Omega^{-1}H_Y\rightarrow 0$
\end{center}
of functors in $\Mod\mathcal{M}$ we have that $\Omega^{-1}H_Y\in \rm Eff(\mathcal{M})^{\bot_{n-1}}=\mathcal{L}_{n-1}(\mathcal{M})$, so $I^1$ that is the injective envelope of $\Omega^{-1}H_Y$ belongs to $\mathcal{L}_{n}(\mathcal{M})$ by the proof of Proposition \ref{P5}. By repeating this argument, in the minimal injective coresolution
\begin{equation}\label{Inj}
0\rightarrow H(Y)\rightarrow I^0 \rightarrow I^1\rightarrow \cdots\rightarrow I^{n} 
\end{equation}
for $H_Y$ in $\Mod\mathcal{M}$ we have $I^0,\ldots,I^{n-1}\in\mathcal{L}_n(\mathcal{M})$ and $\Omega^{-1}H_Y,\ldots,\Omega^{-n+1}H_Y\in\mathcal{L}_1(\mathcal{M})$.
In the last step applying $\Hom(E,-)$ for an arbitrary weakly effaceable functor $E$ to the short exact sequence of functor
\begin{center}
$0\rightarrow \Omega^{-n+1}H_Y\rightarrow I^{n-1}\rightarrow \Omega^{-n}H_Y\rightarrow 0$
\end{center}
in $\Mod\mathcal{M}$ we have the following exact sequence of abelian groups.
\begin{align*}
0&\rightarrow \Hom(E,\Omega^{-n+1}H_Y)\rightarrow\Hom(E,I^{n-1})\rightarrow\Hom(E,\Omega^{-n}H_Y)\\
&\rightarrow\Ext^1(E,\Omega^{-n+1}H_Y)=0.
\end{align*}
Thus $\Hom(\rm Eff(\mathcal{M}),\Omega^{-n}H_Y)=0$, and because $I^n$ is an essential extension of $\Omega^{-n}H_Y$ and $\rm Eff(\mathcal{M})\subseteq \Mod \mathcal{M}$ is a Serre subcategory we have that $\Hom(\rm Eff(\mathcal{M}),I^n)=0$. Therefore by the proof of Proposition \ref{P5} $I^n$ belongs to $\mathcal{L}_n(\mathcal{M})$.
Thus we constructed an injective coresolution \eqref{Inj} for $H_Y$ with $I^0,\ldots,I^n\in \mathcal{L}_1(\mathcal{M})$.
Since the inclusion functor $\mathcal{L}_1(\mathcal{M})\hookrightarrow \Mod\mathcal{M}$ preserve monomorphisms $I^0,\ldots,I^n$ are injective objects in the abelian category $\mathcal{L}_1(\mathcal{M})$.
Thus we have
\begin{center}
$\Ext_{\mathcal{L}_1(\mathcal{M})}^i(H_X,H_Y)\cong\Ext_{\Mod\mathcal{M}}^i(H_X,H_Y)=0,$
\end{center}
 for every $1\leq i\leq n-1$, because representable functors are projective objects in $\Mod\mathcal{M}$.
\end{proof}
\end{proposition}

In the following proposition we prove that the canonical functor $H:\mathcal{M}\rightarrow \mathcal{L}_1(\mathcal{M})$ detect $n$-exact sequences. 

\begin{proposition}\label{P8}
Let $Y:Y^0\rightarrow Y^1\rightarrow \cdots\rightarrow Y^n\rightarrow Y^{n+1}$ be a complex of objects in $\mathcal{M}$ such that
\begin{equation}\label{det}
0\rightarrow H_{Y^0}\rightarrow H_{Y^1}\rightarrow \cdots\rightarrow H_{Y^{n}}\rightarrow H_{Y^{n+1}}\rightarrow 0
\end{equation}
is exact in $\mathcal{L}_1(\mathcal{M})$. Then $Y$ is an admissible $n$-exact sequence in $\mathcal{M}$.
\begin{proof}
Because the essential image of $H:\mathcal{M}\rightarrow \mathcal{L}_1(\mathcal{M})$ is $n$-rigid, by a similar argument like \cite[Proposition 2.2]{JK} for each object $Z\in \mathcal{M}$ we have the following exact sequence of abelian groups.
\begin{equation}
0\rightarrow \Hom(H_Z,H_{Y^0})\rightarrow \Hom(H_Z,H_{Y^1})\rightarrow \cdots\rightarrow \Hom(H_Z,H_{Y^{n}})\rightarrow \Hom(H_Z,H_{Y^{n+1}}) \notag
\end{equation}
Thus by Yoneda's Lemma $Y$ is a left $n$-exact sequence. Dually it is a right $n$-exact sequence, so it is an $n$-exact sequence. We need to show that $Y$ is an admissible $n$-exact sequence. The cokernel of $H_{Y^{n}}\rightarrow H_{Y^{n+1}}$, denoted by $C$, is weakly effaceable. In particular, there exist $X^n\in \mathcal{M}$ and an admissible epimorphism $X^n\twoheadrightarrow Y^{n+1}$ in $\mathcal{M}$, such that $C(Y^{n+1})\rightarrow C(X^n)$ carries the image of $1_{Y^{n+1}}$ to $0$. This means that 
there is a commutative diagram with exact rows in $\mathcal{L}_1(\mathcal{M})$ of the following form for an admissible $n$-exact sequence $X:X^0\rightarrowtail X^1\rightarrow \cdots \rightarrow X^n\twoheadrightarrow Y^{n+1}$ in $\mathcal{M}$. Where the dotted arrows are induced by the factorization property of $n$-kernel.
\begin{center}
\begin{tikzpicture}
%\node (X1) at (-8,1) {$0$};
%\node (X2) at (-6.5,1) {$X^0$};
\node (X3) at (-5,1) {$0$};
\node (X4) at (-3.5,1) {$H_{X^0}$};
\node (X5) at (-1.75,1) {$H_{X^1}$};
\node (X6) at (0,1) {$\cdots$};
\node (X7) at (1.5,1) {$H_{X^n}$};
\node (X8) at (3,1) {$H_{Y^{n+1}}$};
\node (X9) at (4.5,1) {$0$};
\node (X10) at (-5,-0.5) {$0$};
\node (X11) at (-3.5,-0.5) {$H_{Y^0}$};
\node (X12) at (-1.75,-0.5) {$H_{Y^1}$};
\node (X13) at (0,-0.5) {$\cdots$};
\node (X14) at (1.5,-0.5) {$H_{Y^n}$};
\node (X15) at (3,-0.5) {$H_{Y^{n+1}}$};
\node (X16) at (4.5,-0.5) {$0.$};
%\draw [->,thick] (X1) -- (X2) node [midway,left] {};
%\draw [->,thick] (X2) -- (X3) node [midway,left] {};
\draw [->,thick] (X3) -- (X4) node [midway,left] {};
\draw [->,thick] (X4) -- (X5) node [midway,left] {};
\draw [->,thick] (X5) -- (X6) node [midway,left] {};
\draw [->,thick] (X6) -- (X7) node [midway,above] {};
\draw [->,thick] (X7) -- (X8) node [midway,left] {};
\draw [->,thick] (X8) -- (X9) node [midway,above] {};
\draw [->,thick] (X10) -- (X11) node [midway,left] {};
\draw [->,thick] (X11) -- (X12) node [midway,above] {};
\draw [->,thick] (X12) -- (X13) node [midway,above] {};
\draw [->,thick] (X13) -- (X14) node [midway,above] {};
\draw [->,thick] (X14) -- (X15) node [midway,above] {};
\draw [->,thick] (X15) -- (X16) node [midway,above] {};
\draw [->,thick,dotted] (X4) -- (X11) node [midway,above] {};
\draw [->,thick,dotted] (X5) -- (X12) node [midway,above] {};
\draw [->,thick] (X7) -- (X14) node [midway,above] {};
\draw [double,-,thick] (X8) -- (X15) node [midway,above] {};
\end{tikzpicture}
\end{center}
Because the top row is induced by an admissible $n$-exact sequence, by the dual of Obscure axiom (\cite[Proposition 4.11 ]{J}) and the Yoneda's Lemma the bottom row is also induced by an admissible $n$-exact sequence.
\end{proof}
\end{proposition}

\begin{remark}\label{R8}
By the Example \ref{Ex} there are $n$-exact categories that
aren't equivalent to $n$-cluster tilting subcategories. Motivated by the above proposition one can try to prove the following.
\begin{itemize}
\item[$\bigstar$]
Let $\mathcal{M}$ be a small $n$-exact subcategory, is there an exact category $\mathcal{E}$ and an embedding $\mathcal{M}\hookrightarrow \mathcal{E}$ such that the additive closure $\rm add(\mathcal{M})$ is an $n$-cluster tilting subcategory of $\mathcal{E}$? 
\end{itemize}
\end{remark}

By Example \ref{Ex} there are $n$-exact categories that are not $n$-cluster tilting. Every $n$-abelian category has a natural structure of $n$-exact category \cite[Theorem 4.4]{J}. The positive answer to the following question tells that every $n$-exact category can be viewed as a nice subcategory of an $n$-abelian category.

\begin{question}\label{Q}
Let $\mathcal{M}$ be a small $n$-exact category. Is $\rm Eff(\mathcal{M})^{\bot_n}=\mathcal{L}_n(\mathcal{M})$ an
$n$-cluster tilting subcategory of the abelian category $\rm Eff(\mathcal{M})^{\bot_1}=\mathcal{L}_1(\mathcal{M})$?
\end{question}

\begin{remark}
Note that positive answer to Question \ref{Q} complete the following table
in a natural way. Recall that for an additive category $\mathcal{B}$, $\modd \mathcal{B}$ is the full subcategory
of $\rm Mod \mathcal{M}$ consist of all finitely presented functors, and $\rm eff(\mathcal{B})$ is those functors that an epimorphism induces their finite presentation. It is not hard to see that for abelian and $n$-abelian categories $\rm eff(\mathcal{B})=Eff(\mathcal{B})\cap \modd \mathcal{B}$. The first equivalence is called "Auslander's formula". The second equivalence is called "Gabriel-Quillen embedding theorem" (see \cite[Appendix A]{Kl}). And the third equivalence recently was proved in \cite{EN, Kv}.
\begin{center}
\begin{tabular}{|p{2.3in}|p{3.6in}|} \hline 
$\mathcal{A}$ is a small abelian category. & $\dfrac{\rm mod\mathcal{A}}{\rm eff(\mathcal{A})}\simeq \rm eff(\mathcal{A})^{\bot}\simeq\mathcal{A}$\\ \hline 
$\mathcal{E}$ is a small exact category. & $\dfrac{\rm Mod\mathcal{E}}{\rm Eff(\mathcal{E})}\simeq \rm Eff(\mathcal{E})^{\bot}\simeq Lex(\mathcal{E})$, and $\mathcal{E}$ is an extension-closed subcategory of it. \\ \hline 
$\mathcal{M}$ is a small $n$-abelian category. & $\dfrac{\rm mod\mathcal{M}}{\rm eff(\mathcal{M})}\simeq \rm eff(\mathcal{M})^{\bot}$ has an $n$-cluster tilting subcategory equivalent to $\mathcal{M}$.\\ \hline 
$\mathcal{M}$ is a small $n$-exact category. & $\dfrac{\rm Mod\mathcal{M}}{\rm Eff(\mathcal{M})}\simeq \rm Eff(\mathcal{M})^{\bot}\simeq \mathcal{L}_1(\mathcal{M})$ has an $n$-cluster tilting subcategory $\rm Eff(\mathcal{M})^{\bot_n}$ that $\mathcal{M}$ nicely embed in it.\\ \hline 
\end{tabular}
\end{center}
\end{remark}

\section*{acknowledgements}
Special thanks are due to the referee who read this paper carefully and made useful comments and
suggestions that improved the presentation of the paper.

\section*{Funding}
This research was in part supported by a grant from IPM (No. 1400180047).

\end{document}